\documentstyle[amssymb,amsfonts]{amsart}

\def\E{{\hbox{\bf E}}}

 at 10 true pt

\def\beq{\begin{equation}}
\def\eeq{\end{equation}}
\def\barray{\begin{eqnarray*}}
\def\earray{\end{eqnarray*}}
\def\be#1{ \begin{equation}\label{#1} }

\def\bas{\begin{align*}}
\def\eas{\end{align*}}
\def\bi{\begin{itemize}}
\def\ei{\end{itemize}}

\def\Z{{\hbox{\bf Z}}}

\def\eps{\varepsilon}

\def \endprf{\hfill  {\vrule height6pt width6pt depth0pt}\medskip}
\def\emph#1{{\it #1}}
\def\textbf#1{{\bf #1}}

\def\In{{\rm In}}

\def\BZ{{\mathbf Z}}

\def\ep{{\epsilon}}

\theoremstyle{plain}
  \newtheorem{theorem}[subsection]{Theorem}
  \newtheorem{conjecture}[subsection]{Conjecture}
  \newtheorem{problem}[subsection]{Problem}

  \newtheorem{lemma}[subsection]{Lemma}
  \newtheorem{corollary}[subsection]{Corollary}

\theoremstyle{remark}

\theoremstyle{definition}

\include{psfig}

\begin{document}

\title{A structural approach to subset-sum problems   }

\author{  Van Vu}
\address{Department of Mathematics, Rutgers, Piscataway, NJ 08854}
\email{vanvu@@math.rutgers.edu}


\thanks{V. Vu  is supported by  NSF Career Grant 0635606.}

\begin{abstract} We discuss  a structural approach to
subset-sum problems in additive combinatorics. The core of this
approach are Freiman-type structural theorems, many of which will be
presented through the paper. These results have  applications in
various areas, such as number theory, combinatorics and mathematical
physics.

\end{abstract}

\maketitle
\section{Introduction}

Let $A=\{a_{1} , a_{2} , \dots \}$ be a subset of an additive group
$G$ (all groups discussed in this paper will be abelian). Let
$S_{A}$ be the collection of subset sums of $A$

$$
S_A := \{ \sum_{x\in B}x |B\subset A, |B| < \infty \}.
$$

Two related notions that are frequently considered are

$$lA := \{a_{1} + \dots +a_{l} | a_{i} \in A\} $$

$$l^{*}A := \{a_{1} + \dots +a_{l} | a_{i} \in A, {i} \neq {j}\} .$$

We have the trivial relations

$$l^{*} A \subset lA \,\,\,\hbox{\rm and} \,\, \cup _{l} l^{*}A = S_{A}. $$

One can have similar definitions for $A$ being a sequence (repetitions allowed).

{\it Example.}

$A=\{0,1,4 \}$, $G=\Z$, $2A=\{ 0,1,2,4,5,8\}$, $2^{*} A= \{1,4,5
\}$, $S_{A}= \{0,1,4,5 \}$.

$A=\{0,1,4 \}, G=\Z_{5}$, $2A= G$, $2^{*} A= \{0,1,4 \} =S_{A} $.

Now let $A$ be a sequence:

$A=\{1,1, 9 \}$, $G=\Z$, $3A=\{3, 11, 19, 27 \}$,  $3^{*} A= \{11
\}$, $S_{A} =\{1,1,2,9, 10, 11 \}$.

Notice that for a large $l$, $lA$ can be significantly different from $S_{A}$ and $l^{\ast}A$.
In general, it is easier to handle than the later two.

Many basic problems
in additive combinatorics have the following form:

 { \it If $A$ is sufficiently dense in $G$, then $S_{A}$ (or $l^{*} A$ or
$lA$) contains a special element (such as $0$ or a square),  or a
large structure (such as a long arithmetic progression  $G$ itself).
}

The main question is to find the threshold for ``dense''. As
examples, we  present below a few  well-known results/problems in
the area. In the whole paper, we  are going to focus mostly on two
special cases: (1) $G=\Z_{p}$, where $\Z_{p}$ denotes the cyclic
group of residues modulo a large prime $p$; (2) $G = \Z$, the set of
integers.


Following the literature, we say that $A$ is {\it zero-sum-free} if
$0 \notin S_{A}$. Furthermore, $A$ is {\it complete} if $S_{A} =G$
and {\it incomplete} otherwise.  The asymptotic notation is used
under the assumption that $|A| \rightarrow \infty$.

A basic result concerning zero-sum-free sets is the following
theorem of Olson \cite{O2} and Szemer\'edi \cite{SzEH} from the late
1960s, addressing a problem of Erd\H os and Heilbronn \cite{EH}.

\begin{theorem} \label{theorem:szemeredi} (Olson-Szemer\'edi) Let $A$ be a subset of $\Z_{p}$
with cardinality $C\sqrt p $, for a sufficiently large constant $C$.
Then $S_{A}$ contains zero.
\end{theorem}

To see that order $\sqrt p$ is necessary, consider  $A:=\{1,2,
\dots, n \}$, where $n \approx \sqrt {2p}$ is the largest integer
such that $1+ \dots +n < p$.

Concerning completeness, Olson \cite{O1}, proved the following result

\begin{theorem}  \label{theorem:olson}
(Olson)  Any subset  $A$  of $\Z_{p}$ with cardinality at least
$\sqrt {4p-3} +1 $ is complete. \end{theorem}

To see that the bound is close to optimal, take  $A:= \{-m, \dots,
-1,0, 1, \dots, m \}$ where $m$ is the largest integer such that $1+
\dots +m < \lfloor  p/2 \rfloor$.

Another classical result concerning zero sums is that of Erd\H
os-Ginburg-Ziv \cite{EGZ}, again from the 1960s.

\begin{theorem} \label{theorem:EGZ} (Erd\H os-Ginburg-Ziv)  If $A$ is a sequence of $2p-1$
elements in $\Z_{p}$, then $p^{*} A$ contains zero. \end{theorem}

This theorem is sharp by the following example:  $A =\{0^{[p-1]},
1^{[p-1]} \}$ Furthermore, instead of $0$ and $1$, one can use any
two different elements of $\Z_{p}$. (Here and later $x^{[k]}$ means
$x$ appears with multiplicity $k$ in $A$.)

Now we discuss two problems involving  integers. Set $[n]:=\{1, 2,
\dots,  n \}$. An old and popular conjecture concerning subset sums
of integers is Folkman's conjecture, made in 1966 \cite{folkman}.
Folkman's conjecture is a strengthening  of  a conjecture by Erd\H
os \cite{erdfolkman}  about finding a necessary and sufficient
condition for a sequence $A$ such that $S_{A}$ contains all but
finite exception of the positive integers.

\begin{conjecture} \label{conj:folkman} (Folkman's conjecture) The following holds for any
sufficiently large constant $C$. Let $A$ be an strictly increasing
sequence of positive integers with (asymptotic) density at least $C
\sqrt n $ (namely, $|A \cap [n] | \ge C \sqrt n$ for all
sufficiently large $n$). Then $S_{A}$ contains an infinite
arithmetic progression. \end{conjecture}

Cassels \cite{Cass} and Erd\H os \cite{erdfolkman} showed that
density $\sqrt n$ is indeed needed; thus, Folkman's conjecture is
sharp up to the value of $C$.  For more discussion about Folkman's
conjecture and its relation with Erd\H os' conjecture, we refer to
\cite{folkman} and  the monograph \cite{EG} by Erd\H os and Graham.

Finally,  a problem involving a non-linear relation, posed by Erd\H
os in 1986 \cite{erdsquare}.

\begin{problem} \label{conj:erd}  (Erd\H os' square-sum-free problem)
A set $A$ of integers is square-sum-free if $S_A$ does not contain a
square. Find the largest size of a square-sum-free subset of $[n]$.
\end{problem}

Erd\H os observed that one can construct such a square-sum-free
subset of $[n]$ with  at least $\Omega (n^{1/3})$ elements. To see
this, consider $A:= \{q,2q, \dots, kq\}$ with $q $ prime, $(k+1)k <
2q$, $kq  \le n$. Since the sum of all elements of $A$ is less than
$q^{2}$, $S_{A}$ does not contain a square. Erd\H os
\cite{erdsquare} conjectured that the truth is close to this lower
bound.

Problems involving subset sums such as the above (and many others)
have been attacked, with considerable success, using various
techniques: combinatorial, harmonic analysis, algebraic etc. The
reader who is interested in these techniques may want to look at
\cite{Alonhandbook, Sarhandbook, TVbook, Nathansonbook} and the
references therein.

The goal of this paper is to introduce the so-called ``structural
approach'', which has been developed systematically in recent years.
This  approach is  based on the following simple plan

{\it Step I:  Force a structure on $A$.} In this step, one tries to show the following:
If $A$ is relatively dense (close to the conjectured  threshold but not yet there)
and $S_{A}$ does not contain the desired object, then $A$ has a  very special structure.

\noindent Alternately, one can can try to

{\it Step I':  Find a structure in  $S_A$.} If $A$ is relatively
dense (again close to the conjectured threshold but not yet there)
then  $S_{A}$ contains a special structure.

{\it Step II: Completion.} Since $|A|$ is still below the threshold,
we can  add (usually many) new elements to $A$. Using these elements
 together with the existing structure, one can, in most cases,
obtain the desired object in a relatively simple manner.

The success of the method depends on the quality of the information
we can obtain on the structure of $A$ (or $S_A$) in Step I (or I').
In several recent studies, it has turned out that one can frequently
obtain something close to a {\it  complete characterization} of
these sets. Thanks to these results, one is able to make
considerable progresses on  many old problems and also reprove and
strengthen several existing ones (with a better understanding and a
complete classification of the extremal constructions).

The rest of this paper is devoted to the presentation of these
structural theorems and their representative applications.

\section{Freiman's structural theorem}

A corner stone in additive combinatorics is the structural theorem
of Freiman (sometime referred to as Freiman's inverse theorem),
which writes down the  structure of sets with small doubling.

A {\it generalized arithmetic progression} (GAP) of rank $d$ in a
group $G$ is a set of the form

$$\{a_{0}+ a_{1}x_{1} + \dots  +a_{d} x_{d} | M_{i} \le x_{i} \le N_{i}\}, $$

\noindent where $a_{i}$ are elements of $G$ and $M_i\le N_i$ are
integers. It is intuitive to view a GAP $Q$ as the image of the
$d$-dimensional integral  box $B:= \{(x_{1}, \dots, x_{d})| M_{i}
\le x_{i} \le N_{i} \}  $ under the linear map

$$\Phi (x_{1} , \dots, x_{d}) =  a_{0}+ a_{1}x_{1} + \dots  +a_{d} x_{d}. $$

We say that $Q$ is {\it proper} if $\Phi$ is one-to-one. It is easy
to see  that if $Q$ is a proper GAP of rank $d$ and $A$ is a subset
of density $\delta$ of $Q$, then $|2A| \le C(d, \delta) |A|$.
Indeed,

$$|2A| \le |2Q| \le |2B| =2^d |B|=  2^{d} |Q| \le \frac{2^{d}} {\delta} |A| $$
since the volume of a box increases by a factor $2^{d}$ if its sizes are doubled.

Freiman's theorem shows that this is the {\it only} construction of sets with constant doubling.

\begin{theorem} \label{theorem:freiman} (Freiman's theorem)
\cite{freiman} For any positive constant $C$, there are positive constants $d=d(C)$
ad $\delta=\delta(C)$ such that the following holds. Let   $A$ be a finite
 subset of a torsion-free group $G$ such that
$|2A| \le C |A|$. Then there is a proper GAP $Q$ of dimension $d$ such that
$A \subset Q$ and $|A| \ge \delta |Q|$. \end{theorem}

Freiman theorem has been extended recently to the torsion case by
Green and Ruzsa \cite{GR}. \cite[Chapter 5]{TVbook} contains a
detailed discussion of both theorems and related results.

One can use Freiman's theorem iteratively to treat the sumset $lA$ ,
for $l > 2$. For simplicity,  assume that $l=2^{s}$ is a power of
$2$. Thus, the set  $A_s:= lA=2^{s} A$ can be viewed as $2 A_{s-1}$
where $A_{s-1} := 2^{s-1} A$. Using a multi-scale analysis combined
with Fremain's theorem, one can obtain useful structural information
about $lA$ or $A$ itself. For an example of this technique, we refer
to \cite{SzV1} or \cite[Chapter 12]{TVbook}.

The treatment of $l^{\ast} A$ and $S_{A}$ is  more difficult. However,  one can still develop
 structural theorems in these cases. While the content of most theorems in this direction are quite different
 from that of Freiman's, they do bear a similar spirit that somehow
 the most natural construction happens to be (essentially) the only one.

\section{Structure of zero-sum-free sets}

Let $A$ be a zero-sum-free subset of $\Z_{p}$.  We recall the
example following Theorem \ref{theorem:szemeredi}. Let   $A:=\{1,2,
\dots,  n  \}$. If $1+ \dots +n < p$, then obviously  $S_{A}$ does
not contain $0$. This shows that a zero-sum-free set can have close
to $\sqrt {2p}$ elements. In \cite{SzV1}, Szemer\'edi and Vu showed
that {\it having elements with small sum is essentially the only
reason for a set to be zero-sum-free}.  More quantitative versions
of this statement were worked out in \cite{NSV} and \cite{NV1}. For
example, we have \cite[Theorem 2.2]{NV1}

\begin{theorem} \label{theorem:NV1} After a proper dilation (by some non-zero element), any
zero-sum-free subset $A$ of $\Z_{p}$  has the form

$$A = A' \cup A^{''} $$

\noindent where the elements of $A'$ (viewed as integers between 0
and $p-1$) are small, $\sum_{x \in A'} x < p$, and $A^{''}$ is
negligible, $|A^{''}| \le p^{6/13+o(1)}$. \end{theorem}

One can perhaps improve the constant $6/13$ by tightening the analysis in
\cite{NV1}. It is not clear, however, what would be the best constant here.
In most applications, it suffices to have any constant strictly less than $1/2$.

The dilation is necessary. Notice that if $A$ is zero-sum-free
(incomplete), then the set  $A_{x}:= \{xa| a \in A \}$ is also
zero-sum-free  (incomplete) for any $0\neq x  \in \Z_{p}$.

We can also prove similar results for $lA$ and $l^{*}A$, and for $A$
being a sequence (see \cite{NV1} for details). In the rest of this
section, we present few applications of these results.

\subsection{The size of the largest zero-sum-free set in $\Z_{p}$}

Let $m_{p}$ denote the size of the largest zero-sum-free set in
$\Z_{p}$. The problem of determining $m_{p}$ was posed by Erd\H os
and Heilbronn \cite{EH} and has a long history. Szemer\'edi proved
that $m_{p} \le C\sqrt p$, for some sufficiently large $C$
independent of $p$ \cite{SzEH}. Olson showed that $C=2$ suffices
\cite{O2}.
 Much later, Hamidoune and Z\'emor \cite{HZ} showed that
$m_{p} \le \sqrt {2p} + 5 \log p$, which is asymptotically sharp. Using an earlier version of Theorem \ref{theorem:NV1},
 Szemer\'edi, Nguyen and Vu \cite{NSV}  recently obtained the exact value of $m_{p}$.

\begin{theorem}  \label{theorem:SNV1} Let $n_{p}$ be the largest integer
so that $1+\dots+ (n_{p}-1) < p$.
\begin{itemize}

\item If $p \neq \frac{n_{p}(n_{p} +1)}{2} -1$,  then $m_{p}=n_{p}-1$.

 \item If $p = \frac{n_{p}(n_{p}+1)}{2} -1$, then $m_{p} =n_{p}$. Furthermore, up
to a dilation, the only zero-sum-free set with $n_{p}$ elements is
$\{-2,1,3,4, \dots, n_{p} \}$.

\end{itemize}
 \end{theorem}

The same result was obtained by Deshouillers and Prakash (personal
communication by Deshouillers) at about the same time.

\subsection{The structure of relatively large zero-sum-free sets}

Let us now consider the structure of zero-sum-free sets of size close to $\sqrt {2p}$.
Let $\|x\|$ denote the integer norm of $x$. In \cite{Des1}, Deshouillers proved

\begin{theorem} \label{theorem:Des1} \label{theorem:Des1} Let $A$ be a
zero-sum-free subset of $Z_{p}$ of size at least $ \sqrt {p}$. Then
(after a proper dilation)

$$\sum_{x\in A, x <p/2} \|x/p\|  \le  1 + O(p^{-1/4} \log p) $$

$$\sum_{x\in A, x >p/2} \|x/p\|  \le   O(p^{-1/4} \log p). $$ \end{theorem}

Deshouillers showed (by a construction) that the error term
$p^{-1/4}$ cannot be replace by
 $o(p^{-1/2})$.  Using  an earlier  version of Theorem \ref{theorem:NV1},
 Nguyen, Szemer\'edi and Vu \cite{NSV} improved Theorem \ref{theorem:Des1}
  to  obtain the best possible error term
  $O(p^{-1/2})$, under a stronger assumption on the size of $|A|$.

 \begin{theorem} \label{theorem:Des1} \label{theorem:Des1-1}  Let $A$ be a zero-sum-free subset of $Z_{p}$ of size at least $ .99 \sqrt {2p}$. Then (after a proper dilation)

$$\sum_{x\in A, x <p/2} \|x/p\|  \le  1 + O(p^{-1/2}) $$

$$\sum_{x\in A, x >p/2} \|x/p\|  \le   O(p^{-1/2}). $$ \end{theorem}

The constant $.99$ is, of course, ad-hoc and can be improved by redoing the analysis carefully.
 On the other hand,  it is not clear what  the best assumption on $|A|$ should be.

\subsection{Erd\H os-Ginburg-Ziv revisited}

Using a version of Theorem \ref{theorem:NV1} for sequences, Nguyen
and Vu \cite{NV1} obtained the following characterization for a
sequence of size slightly more than $p$ and does not contain a
subsequence of $p$ elements summing up to $0$.

\begin{corollary}\label{cor:EGZ1} \cite[Theorem 6.2]{NV1} Let $\eps$
be an arbitrary positive constant.
 Assume that $A$ is a $p$-zero-sum-free sequence
and $p+ p^{12/13+\epsilon}\le |A| \le 2p-2$. Then $A$ contains two
elements $a$ and $b$ with multiplicities $m_{a},m_{b}$ satisfying
$m_a+m_b \ge 2(|A|-p- p^{12/13 +\epsilon})$.
\end{corollary}

The interesting point here is that the structure kicks in very soon,
when $A$ has just slightly more than $p$ elements. Few years ago,
Gao, Panigrahi, and Thangdurai \cite{GPT} proved a similar statement
under the stronger assumption that $|A| \ge 3p/2$.

 It is easy to deduce Erd\H os-Ginburg-Ziv theorem from Corollary \ref{cor:EGZ1}, together with a complete
 characterization of the extremal sets. The reader may want to consider
 as an exercise or check \cite{NV1} for  details.

\section{ Incomplete Sets}

Now we turn our attention to incomplete sets, namely sets $A$ where
$S_{A} \neq \Z_{p}$. The situation here is very similar to that with
zero-sum-free sets.  Szemer\'edi and Vu
 \cite{SzV1}  showed that {\it having elements with small sum is essentially
 the only reason for a set to be incomplete}.
 More quantitative  versions of this statement were worked out in \cite{NSV} and \cite{NV1}.
 For example, in \cite{NV1}, the following analogue of Theorem \ref{theorem:NV1} was proved

\begin{theorem}  \label{theorem:NV2}  After a proper dilation (by some non-zero element), any
incomplete subset  $A$ of $Z_{p}$  has the form

$$A = A' \cup A^{''} $$

where the elements of $A'$ are small
(in the integer norm), $\sum_{x \in A'}  \| x /p \| < 1$
and $A^{''}$ is negligible, $|A^{''}| \le p^{6/13+o(1)} $. \end{theorem}

The reader can find similar  results for $lA$ and $l^{*}A$ and for
$A$ being a sequence in \cite{NV1}. We next discuss some
applications of these results.

\subsection{The structure of relatively large incomplete sets}

Theorem \ref{theorem:NV2} enables us to prove results similar to
those in the last section for incomplete sets. The problem of
determining the largest size of an incomplete set in $\Z_{p}$ was
first considered by Erd\H os and Heilbronn \cite{EH} and essentially
solved by Olson (Theorem \ref{theorem:olson}).
 da Silva and Hamidoune  \cite{daH} tightened the bound to
$\sqrt {4p-7} + 1$, which is best possible. We are not going to go
into these results here,
 but note that one can perhaps obtain a new proof (with a characterization of the extremal sets)
  using Theorem \ref{theorem:NV2}.

Concerning the structure of relatively large incomplete sets,
Deshouillers and Freiman  \cite{DF} proved

\begin{theorem} \label{theorem:Des1} \label{theorem:DF} Let $A$ be an incomplete subset of $\Z_{p}$
of size at least $ \sqrt {2p}$. Then (after a proper dilation)

$$\sum_{x\in A} \|x/p\|  \le  1 + O(p^{-1/4} \log p) .$$
 \end{theorem}

They conjectured that the error term may be replaced by $O(\sqrt p)$, which would be best possible due to a later construction of Deshouillers \cite{Des2}.

Using Theorem \ref{theorem:NV2}, Nguyen and Vu \cite{NV1} confirmed
this conjecture, provided that $A$ is sufficiently close to $2\sqrt
p$.

\begin{theorem} \label{theorem:Des1} \label{theorem:DF1} Let $A$ be an incomplete subset of $Z_{p}$ of size at
least $ 1.99 \sqrt p$. Then (after a proper dilation)

$$\sum_{x\in A} \|x/p\|  \le  1 + O(\sqrt p) .$$
 \end{theorem}

Similar to the constant $.99$ (in the previous section), the
constant $1.99$ is ad-hoc and can be improved by redoing the
analysis carefully. On the other hand,  it is not clear what  the
best assumption on $|A|$ is.

\subsection{The structure of incomplete sequences}

Let us now discuss (rather briefly) the situation with sequences.
The main difference between sets and sequences is that a sequence
can have elements with high multiplicities. It has turned out that
when the maximum multiplicity of incomplete sequence $A$ is
determined, one can obtain  strong structural information about $A$.

Let $1\le m \le p$ be a positive integer and  $A$ be an incomplete
sequence of $\Z_p$ with maximum multiplicity $ m$. Trying to make
$A$ as large as possible, we come up with the following example,

$$B_m=\{-n^{[k]}, (n-1)^{[m]}, \dots,-1^{[m]},0^{[m]},1^{[m]},\dots,(n-1)^{[m]}, n^{[k]}\}$$

where $1 \le k \le   m$ and $n$ are the unique integers satisfying

$$2m(1+2+\dots+n-1)+2kn<p\le 2m(1+2+\dots+n-1)+2(k+1)n. $$

It is clear that any subsequence of $B_m$ is incomplete and has
multiplicity at most $m$. In \cite{NV1}, we proved that any
incomplete sequence $A$ with maximum multiplicity $m$ and
cardinality close to $|B_m|$ is essentially a subset of $B_m$.

\begin{theorem}\label{theorem:NV3} Let $6/13< \alpha< 1/2$
be a constant. Assume that $A$ is an incomplete sequence of $\Z_p$
with maximum multiplicity $m$ and cardinality $|A|=
|B_m|-O((pm)^{\alpha})$. Then after a proper dilation, we can have

$$ A = A' \cup A^{''} $$

where $A' \subset B_m$ and $|A^{''}|=O((pm)^{(\alpha+1/2)/2})$.
\end{theorem}

\subsection{Counting problems}

Sometime one would like to count the number of sets which forbid
certain additive configurations. A well-known example of problems of
this type is  the Cameron-Erd\H os problem \cite{EC}, which asked
for the number of subsets of $[n]= \{1,2, \dots, n \}$ which does
not contain three different elements $x,y,z$ such that $x+y-z=0$.
Cameron an Erd\H os noticed that any set of odd numbers has this
property. Thus, in $ [n]$ there are  at least $\Omega (2^{n/2})$
subsets with the required property. They conjectured that $2^{n/2}$
is the right order of magnitude. There were several partial results
\cite{AlonCE, Calkin, ErdGran} before Green settled the conjecture
\cite{green}.

Using structural theorems such as Theorem \ref{theorem:NV1}, one can
obtain results of similar spirit for the number of zero-sum-free or
incomplete  sets and sequences. For example, using an earlier
version of Theorem \ref{theorem:NV1} and standard facts from the
theory of partitions \cite{A}, Szemer\'edi and Vu \cite{SzV1} proved

 \begin{corollary}\label{cor:incomplet:partition:set}
 The number of incomplete subsets of $\Z_{p}$ is
 $\exp((\sqrt{\frac{2}{3}}\pi + o(1))\sqrt{p})$.
 \end{corollary}

Using Theorem \ref{theorem:NV3}, one obtains the following
generalizations \cite{NV1}.

 \begin{corollary}\label{theorem:incomplete:partition:m}
 The number of incomplete sequences $A$ with highest multiplicity $m$ in $\Z_p$ is
 $\exp((\sqrt{(1-\frac{1}{m+1})\frac{4}{3}}\pi + o(1))\sqrt{p}).$
\end{corollary}

It is an interesting question to determine the error term $o(1)$.

\section{Incomplete sets in a general abelian group}

Let us now consider the problem of finding the largest size of an
incomplete set in a general abelian group $G$, which we denote by
$\In(G)$ in the rest of this section. The situation with a general
group $G$ is quite different from that with $\Z_{p}$, due to the
existence of non-trivial subgroups. It is clear that any such
subgroup is incomplete.  Thus, $\In(G) \ge h$, where $h$ is the
largest non-trivial divisor of $|G|$. The intuition behind the
discussion in this section is that {\it a large incomplete set
should be essentially contained in a proper subgroup}.

In 1975, Diderrich \cite{D1} conjectured that if $|G|=ph$, where $p
\ge 3$ is the smallest prime divisor of $|G|$ and $h$ is composite,
then $c(G) = h+p-2$. (The cases where $p=2$ or $h$ is a prime is
simpler and were treated earlier, some by Diderrich himself
\cite{D1, MW, DM}.)  Didderich's conjecture was settled by Gao and
Hamidoune in 1999 \cite{GH}.

The following simple fact explains the appearance of the  term $p-2$.

{\bf Fact.} If $S_{A \cap H} =H$ for some maximal subgroup $H$ of (prime)  index $q$, then
$|A| \le |H| + q-2$.

  To verify this fact, notice that $A/H$ is a sequence in the  group $\Z_{q}$. It is easy to
  show (exercise)
that  if
 $B$ is a sequence of  $q-1$ non-zero elements in $\Z_{q}$, then $S_{B} \cup \{0\}= \Z_{q}$.

We say that  subset $A$ of $G$ is {\it sub-complete} if there is a
subgroup $H$ of prime index such that $S_{A\cap H} = H$.

Once we know that an incomplete set  $A$ is sub-complete, we can
write down its structure completely.  There is a subgroup $H$ with
prime index $q$ such that $|A \backslash H | \le q-2$, and the
sequence $A/H$ is incomplete in $\Z_q$. (The structure of such a
sequence was discussed in the previous section.) It is natural to
pose the following

\begin{problem} Find the threshold for sub-completeness. \end{problem}

Recently, Gao, Hamidoune, Llad\'o and Serra  \cite{GHLS}  showed
(under some weak assumption) that  any subset of at least
$\frac{p}{p+2}  h +p$ elements is  sub-complete. Furthermore, one
can choose $H$ to have index $p$, where $p$ is the smallest prime
divisor of $|G|$.  Vu \cite{Vudir} showed (again under some weak
assumption)that $\frac{5}{6}h$ is sufficient to guarantee
sub-completeness. It is not clear, however, that what the sharp
bound is.

The situation is much better if we assume that $|G|$ is sufficiently
composite. In particular, if  the product of the two smallest prime
divisors of $|G|$ is significantly smaller than $\sqrt {|G|}$, then
one can determine the sharp threshold for sub-completeness.

\begin{theorem} \label{theorem:dir2} \cite{Vudir} For any positive constant  $\delta$
there is  a positive constant $D(\delta)$ such that the following
holds. Assume that $|G| =p_1 \dots p_t$ , where $t\ge 3$ and  $p_1
\le p_2 \dots \le p_t$ are  primes such that $p_1 p_2\le
\frac{1}{D(\delta)} \sqrt {|G|/ \log |G|}$. Then any incomplete
subset  $A$ of $G$ with cardinality at least $(1+
\delta)\frac{n}{p_1 p_2}$ is subcomplete. Furthermore, the lower
bound $(1+ \delta)\frac{n}{p_1 p_2}$ cannot be replaced by
$\frac{n}{p_1 p_2} + n^{1/4-\alpha}$, for any constant $\alpha$.
\end{theorem}


\section{Structures in $S_{A}$}

As mentioned in the introduction, an alternative way to implement our plan is to find a structure in $S_{A}$ rather than in $A$ (Step I').  A well-known result concerning the structure of $S_{A}$ is the following theorem, proved by
Freiman \cite{freiman2} and S\'ark\"ozy \cite{sar1} independently.

\begin{theorem}  \label{theorem:FS} There are positive constants $C$ and $c$ such that the  following holds for all sufficiently large $n$. Let $A$ be a subset of $[n]:=\{1, \dots, n \}$ with at least
$C \sqrt {n \log n}$ elements. Then $S_{A}$ contains an arithmetic progression of length $c |A|^{2} $.
\end{theorem}

It is clear that the bound on the length of the arithmetic
progression (AP) is optimal, as one can take $A$ to be an interval.
The lower bound on $|A|$, however, can be improved to $C \sqrt n $,
as showed by Szemer\'edi and Vu \cite{SzV2}.

\begin{theorem}  \label{theorem:SzemVu2} There are positive constants $C$ and $c$ such that the  following holds for all sufficiently large $n$. Let $A$ be a subset of $[n]:=\{1, \dots, n \}$ with at least  $C \sqrt {n }$ elements. Then $S_{A}$ contains an arithmetic progression of length $c |A|^{2} $.
\end{theorem}

The assumption $|A| \ge C \sqrt n$ is optimal, up to the value of $C$, as one can construct
a set $A \subset [n]$ of $\epsilon \sqrt n$ elements, for some small constant $\eps$, such that $S_{A}$ does not contain any arithmetic progression of length larger than $n^{3/4} $ (see \cite{SzV2} or
\cite[Section 3.4]{SzV3}).

Theorem \ref{theorem:SzemVu2} can be extended considerably.
Szemer\'edi and Vu \cite{SzV3} showed that for any set $A \subset
[n]$ and any integer $l$ such that $l ^{d} |A| \ge n$ for some
constant $d$, the sumset $l^{\ast} A$ contains a large proper
generalized arithmetic progression (GAP). The parameters of this GAP
is optimal, up to a constant factor (see \cite[Section 3]{SzV3} for
more details).

\begin{theorem} \label{theorem:SzemVu3}  \cite[Theorem 7.1]{SzV3} For any fixed positive integer $d$
there are positive constants $C$ and $c$ depending on $d$
 such that the following holds. For any positive integers $n $
and $l$ and any set  $A \subset [n]$ satisfying $l^{d} |A| \ge C
n$,  $l^{\ast}A$ contains a proper GAP of rank $d'$ and volume at least
$cl^{d'} |A|$, for some integer $1 \le d' \le d$.
\end{theorem}

There are variants of Theorem \ref{theorem:SzemVu3} for finite
fields, and also for sums of different sets (see \cite[Section
5]{SzV3} and \cite[Section 10]{SzV3}). In the following subsections,
we discuss few applications of Theorems \ref{theorem:SzemVu2} and
\ref{theorem:SzemVu3}.

\subsection{Folkman conjectures on infinite arithmetic progressions}

Let us recall to the conjecture of Folkman, mentioned in the introduction.

\begin{conjecture} \label{conj:folkman} (Folkman's conjecture)
The following holds for any  sufficiently large constant $C$. Let
$A$ be an strictly increasing sequence of positive integers with
(asymptotic) density $A(n)$ at least $C \sqrt {n} $ (namely $A(n):=
|A \cap [n] | \ge C \sqrt n$ for all sufficiently large $n$). Then
$S_{A}$ contains an infinite arithmetic progression.
\end{conjecture}

 Folkman \cite{folkman} showed that the conjecture holds
under a stronger assumption  that $A(n) \ge n^{ 1/2 +\ep}$, where
$\ep$ is an arbitrarily small positive constant. (An earlier result
of Erd\H os \cite{erdfolkman} on a closely related problem can
perhaps be adapted to give a weaker bound $n^{(\sqrt 5-1)/2}$.)
 Hegyv\'ari \cite{Heg} and \L uczak and Schoen \cite{LS},
 independently, reduced the density $n^{1+\ep}$ to $C \sqrt{n \log n}$, using
Theorem \ref{theorem:FS}.

Using the stronger Theorem \ref{theorem:SzemVu2}, together with some
additional  arguments, Szemer\'edi and Vu \cite{SzV2} proved the
full conjecture.

\begin{theorem} \label{theorem:folkman} Conjecture \ref{conj:folkman} holds. \end{theorem}

In the same paper \cite{folkman}, Folkman also made a related
conjecture for increasing, but not strictly increasing sequences.
Let $A(n)$ now be the number of elements of $A$ (counting
multiplicities) at most $n$.

\begin{conjecture} \label{conj:folkman2} (Folkman's second
conjecture) The following holds for any sufficiently large constant
$C$. Let $A$ be an  increasing sequence of positive integers with
such that  $A(n) \ge Cn$ for all sufficiently large $n$. Then
$S_{A}$ contains an infinite arithmetic progression.
\end{conjecture}

Despite  the huge change  from $\sqrt n$ to $n$ in the density
bound, this conjecture is also  sharp \cite{folkman}, and (for some
time) appeared  more subtle  than the first one (see a discussion in
\cite[Chapter 6]{EG}). Folkman \cite{folkman} proved the conjecture
under the stronger assumption that $A(n) \ge n^{1+\eps}$. It does
not seem  that  one can obtain the analogue of Hegyv\'ari and \L
uczak-Schoen results due to the lack of a "sequence" variant of
Theorem \ref{theorem:FS}. However, the method in \cite{SzV3} is
sufficiently robust  to enable one to obtain such a  variant for the
stronger  Theorem \ref{theorem:SzemVu2}. With the help of this
result, one can settle Conjecture \ref{conj:folkman2}

\begin{theorem} \label{theorem:folkman2}   \cite[Section 6]{SzV3}  Conjecture \ref{conj:folkman2} holds. \end{theorem}

The strategy for the proofs of Theorems \ref{theorem:folkman} and \ref{theorem:folkman2} is the following. We first find a sufficient condition for a sequence $A$ such that $S_{A}$ contains an infinite
AP.

We say that an infinite  sequence $A$ admits a {\it good} partition if it can
be partitioned into two subsequences $A'$ and $A^{''}$ with the
following two properties

\begin{itemize}

\item There is a number $d$ such that   $S_{A'}$ contains an
arbitrary long arithmetic progression with difference $d$.

\item Let $A^{''} = b_1 \le b_2 \le b_3 \le \dots$. For any number
$K$, there is an index $i(K)$ such that $\sum_{j=1}^{i-1} b_j \ge
b_i +K$ for all $i \ge i(K)$.

\end{itemize}

\begin{lemma} \label{lemma:goodpartition} If a  sequence $A$  admits a
good partition then $S_{A}$ contains an infinite AP . \end{lemma}

 The second assumption is easy to satisfied given that $A$ has proper density.  Thus, the key is the first assumption.
The main feature here is that in this assumption,  we only need to guarantee the existence of long (but finite) APs. So, Theorem \ref{theorem:SzemVu2}
and its variants can be used with full power to achieve this goal.

\subsection{Erd\H os conjecture on square-sum-free sets}

In this section, we return to Erd\H os conjecture on square-sum-free
sets, mentioned in the introduction. Let $SF (n)$ denote the size of
the largest subset $A$ of $[n]$ such that $S_{A}$ does not contain a
square (or $A$ is square-sum-free). Erd\H os \cite{erdsquare}
observed that $SF(n) =\Omega (n^{1/3})$ and conjectured that the
truth is close to this lower bound. Since then, there have been
several attempts on his conjecture.
  Alon \cite{Alon} proved  that $SF(n) =O( \frac{n}{\log n}).$
 In \cite{Lipkin}  Lipkin improved the bound to
$SF(n)=O( n^{3/4+\varepsilon}).$ Later,  Alon and  Freiman
\cite{AlonFreiman}  obtained another improvement $SF(n) =O(
n^{2/3+\varepsilon}). $ About  fifteen years ago,   S\'ark\"ozy
\cite{Sarkozy} showed $SF(n)= O( \sqrt {n\log n}).$

Let us now address the problem from our  structural approach point
of view. Theorem \ref{theorem:SzemVu2} is no longer useful, as we
are dealing with sets of size around $n^{1/3} $, way below the lower
bound $\sqrt n $ required in this theorem. Fortunately, we have a
more general result, Theorem \ref{theorem:SzemVu3}, which enables us
to find structures in $S_{A}$ for any set of size $n^{\delta}$, for
any constant $\delta$. In particular, we can deduce from this
theorem the following corollary.

\begin{corollary} There are positive constants $C$ and $c$ such that the following holds for all sufficiently large $n$.
Let $A $ be a subset of $[n]$ with cardinality at least $C n^{1/3}$.
Then $S_{A}$ contains either an AP of length $c|A|^{2} $ or a proper
GAP of rank $2$ and volume $c |A|^{3}$.

\end{corollary}

Combining this corollary with some number theoretic arguments,  Nguyen and the author  \cite{NV2} can get
close to the conjectured bound.

\begin{theorem}\label{theorem:square:1} There is a constant $C$ such that $SF(n) \le n^{1/3} \log^C n$. \end{theorem}

We strongly believe that the  $\log$ term can be removed. Details
will appear elsewhere.

\section {Inverse Littlewood-Offord theorems and  random matrices}

In this final section, we discuss a problem with a slightly different nature.
Let $A$ be a sequence of non-zero integers. Now we are going to view $S_{A}$ as a multi-set of $2^{n}$ elements. We denote by
$M_{A}$ be the largest multiplicity in $S_{A}$. For example, if $A=\{1, \dots, 1\}$ , then $M_{A} =
{n \choose \lfloor n/2 \rfloor } = \Theta (2^{n} / \sqrt n).$

The problem of bounding $M_{A}$ originated from Littlewood and
Offord's   work on random polynomials \cite{LO}. In particular, they
proved that  $M_A= O(2^{n} \log n / \sqrt n). $ The $\log $ term was
removed by Erd\H os \cite{erdLO}, who obtained a sharp bound for
$M_A$. Many extensions of this result were obtained by various
researcher: Erd\H os-Moser \cite{EM}, S\'ark\"ozy-Szemer\'edi
\cite{SS},
  Katona \cite{Kat}  Kleitman \cite{Kle},   Hal\'asz \cite{Hal}, Griggs et. al.
  \cite{GLOS}, Frankl-F\"uredi \cite{FF}, Stanley \cite{Stan} etc.
  Among others, it was showed that the bound on $M_{A}$ keeps
  improving, if one forbids more and more additive structures in
  $A$.
  For example,  Erd\H os and Moser
  \cite{EM} showed that if  the  elements of $A$ are different (i.e., $A$ is a set),
  then $M_{A} = O(2^{n} \log n/ n^{3/2} )$. In general, the following can be deduced from results of
   \cite{Hal} (see also \cite[Problem 7.2.8]{TVbook})

  \begin{theorem} \label{theorem:halasz} For any fixed integer $k$ there is a constant $C$ such that the following holds. Let
   $A=\{a_1, \dots, a_n \}$ and  $R_k$ be the number of roots of the equation $$\eps_1 a_{i_1} + \dots +\eps_{2k}
  a_{i_{2k}}=0$$ with $\eps_i = \pm 1$ and $i_1, \dots, i_{2k} \in
  [n]$. Then $M_A \le C n^{-2k-1/2} R_k $.
\end{theorem}

In \cite{TVinverse}, Tao and Vu introduced the notion of Inverse
Littlewood-Offord theorems. The intuition  here is that if  $M_{A}$
is large (of order $2^{n}/n^{C}$ for any constant $C$, say), then
$A$ should have  a very strong structure.

The most general example we found with large $M_{A}$ is the
following. Let $Q$ be a proper GAP of constant rank $d$ and volume
$V$.  If $A$ is a subset of $Q$, it is easy to show that $M_{A} =
\Omega (\frac{1}{n^{d/2} V })$ (in order to see this, view the
elements of $S_{A} $ as random sums $\sum_{i=1} ^{n} \xi_i a_{i}$
where $a_{i}$ are elements of $A$ and $\xi_{i}$ are iid random
variables taking values $0$ and $1$ with probability $1/2$). Thus,
if the volume of $Q$ is small, then $M_{A}$ is large.

 In
\cite{TVinverse}, Tao and Vu proved the inverse statement, asserting
that  {\it having $A$ as a subset of a small GAP is essentially the
only way to guarantee make $M_{A}$ large}.

\begin{theorem}  \label{theorem:TVinverse} \cite{TVinverse}
For any constant $C$ and $\epsilon$  there are constants $B$ and $d$
such that the following holds. Let $A$ be a sequence of $n$ elements
in a torsion-free group $G$. If $M_{A} \ge 2^{n}/n^{C}$ for some
constant $C$, then all but at most  $n^{1-\epsilon}$ elements of $A$
is contained in a proper GAP $Q$ of rank $d$ and cardinality
$n^{B}$.
\end{theorem}

In a more recent paper \cite{TVinverse08},  the same authors
obtained a (near) optimal relationship between the parameters $C,
\epsilon, d$ and $B$. As a corollary, one can deduce (asymptotic
versions of) many earlier results, such as Theorem
\ref{theorem:halasz}.
 (In spirit, this process is somewhat similar to the process
 of using Theorem \ref{theorem:NV1} to reprove, say,  Erd\H os-Ginburg-Ziv
 theorem.)

We would like to conclude this survey with a rather unexpected
application. Let us leave combinatorial number theory and jump to
the (fairy remote) area of mathematical physics. In the 1950s,
Wigner observed and proved his famous semi-circle law concerning the
limiting distribution of eigenvalues in a {\it symmetric} random
matrix \cite{Wig}. A brother of this law, the so-called circular law
for {\it non-symmetric} random matrices,  has been conjectured, but
remains open since that  time.

\begin{conjecture} (Circular Law Conjecture) Let $\xi$ be a random variable with mean 0 and variance 1 and
$M_{n}$ be the random matrix whose entries are iid copies of $\xi$.
Then the limiting distribution of the  eigenvalues of
$\frac{1}{\sqrt n}M_{n}$ converges to the uniform distribution on
the unit disk. \end{conjecture}

Girko \cite{Girko} and Bai \cite{Bai} obtained important partial
results concerning this conjecture. These results and many related
results are carefully discussed in the book \cite{BS}. There has
been a series of rapid developments  recently by G\"otze-Tikhomirov
\cite{GT, GT2},  Pan-Zhou \cite{PZ}, and Tao-Vu \cite{TVcir}. In
particular,  Tao and Vu \cite{TVcir} confirmed the conjecture under
the slightly stronger assumption that the $(2+\eta)$-moment of $\xi$
is bounded, for any $\eta >0$.

\begin{theorem} \label{theorem:CL}
The Circular Law holds (with strong convergence) under an additional
assumption that
 $$\E
(|\xi|^{2+\eta}) < \infty$$ for some fixed $\eta >0$.
\end{theorem}

The key element of this proof is   a variant  of Theorem
\ref{theorem:TVinverse}, which enables us to count the number of
sequences $A$ with bounded elements such that $M_A$ (more precisely
a continuous version of it) is large. For  details, we refer to
\cite{TVcir}.


\begin{thebibliography}{99}


\bibitem{A} G. E. Andrews, { The theory of partitions.} Cambridge
university press, 1998.
\bibitem{AlonCE} N. Alon, Independent sets in regular graphs and
sum-free subsets of abelian groups, {\it Israel Journal Math.} 73
(1991), 247-256.


\bibitem{Alonhandbook}  N. Alon,   Combinatorial Nullstellensatz,
Recent trend in combinatorics (M\'atrah\'aza, 1995), {\it Combin.
Probab. Comput. } 8 (1999), 7-29.


\bibitem{Alon} {N. Alon}, { Subset sums.} {\it Journal of Number Theory,} 27 (1987),
196-205.

\bibitem{AlonFreiman} {N. Alon } and {G. Freiman}, { On sums of subsets of a set of
integers}, {\it Combinatorica}, 8 (1988), 297-306.



 \bibitem{Bai} Z. D. Bai, Circular law, {\it  Ann. Probab. } 25  (1997),  no. 1,
494--529.

\bibitem{BS}
Z. D. Bai and J. Silverstein, Spectral analysis of large dimensional
random matrices, Mathematics Monograph Series \textbf{2}, Science
Press, Beijing 2006.



\bibitem{BD} {A. Bialostocki and P. Dierker}, { On Erd\H{o}s-GinzburgiZiv theorem and the
Ramsey number for stars and matchings}, {\it Discrete Mathematics},
110, (1992), 1-8.

\bibitem{Calkin} N. Calkin, On the number of sum-free sets, {\it
Bull. London Math. Soc.} 22 (1990), 141-144.


\bibitem {Cass}  J.W.S Cassels,  On the representation of integers
as the sums of distinct summands taken from a fixed set, {\it Acta
Sci. Math. Szeged}  21 1960 111--124.




\bibitem{EC}  P. Cameron and P. Erd\H os, On the number of sets of
integers with various properties, Number Thoery (Banff, AB 1988),
61-79, de Gruyter, Berlin, 1990.

\bibitem{daH} D. da Silva and Y. O. Hamidoune,
 Cyclic spaces for Grassmann derivatives and additive theory,
{\it Bull. London Math. Soc.}  26 (1994), no. 2, 140--146.



\bibitem{D1} G. T. Diderrich,
An addition theorem for abelian groups of order $pq$, {\it J. Number
Theory} 7 (1975), 33--48.

\bibitem{DM} G. T. Diderrich and H. B.
 Mann,  Combinatorial problems in finite Abelian groups,
 {\it  Survey of combinatorial theory } (Proc. Internat. Sympos., Colorado State Univ.,
 Fort Collins, Colo., 1971), pp. 95--100. North-Holland, Amsterdam,
 1973.

\bibitem{Des1} J-M. Deshouilers, Quand seule la sous-somme vide est nulle modulo $p$. (French)
[When only the empty subsum is zero modulo $p$]  {\it J. ThŽor.
Nombres Bordeaux } 19  (2007),  no. 1, 71--79.


\bibitem{Des2} { J-M Deshouillers}, {Lower bound concerning subset sum
wich do not cover all the residues modulo $p$}, {\it Hardy-
Ramanujan Journal}, Vol. 28(2005) 30-34.


\bibitem{DF}
{ J-M Deshouillers} and {G. Freiman}, { When subset-sums do not
cover all the residues modulo $p$}, {\it Journal of Number Theory }
104(2004) 255-262.


\bibitem{erdLO}
P. Erd\H {o}s, On a lemma of Littlewood and Offord, \emph{Bull.
Amer. Math. Soc.} \textbf{51} (1945), 898--902.


\bibitem{erdsquare} P. Erd\H os, Some problems and results on combinatorial number theory,
Proceeding of the first China conference in Combinatorics, 1986.


\bibitem {erdfolkman} P. Erd\H os, On the representation of large
interges as sums of distinct summands taken from a fixed set, {\it
Acta. Arith.} 7 (1962), 345-354.

\bibitem {EG} P. Erd\H os and R. Graham,  Old and new problems and results in
combinatorial number theory. Monographies de L'Enseignement
Mathématique [Monographs of L'Enseignement Mathématique], 28.
Université de Genève, L'Enseignement Mathématique, Geneva, 1980.

\bibitem{ErdGran} P. erd\H os and A. Granville, {\it Unpublished.}



\bibitem{EH} P. Erd\H os and H. Heilbronn, On the addition of residue classes ${\rm mod}
p$, {\it Acta Arith.}  9 1964 149--159.

\bibitem{EM} P. Erd\H os and L. Moser,
P. Erd\H{o}s, Extremal problems in number theory. 1965 Proc. Sympos.
Pure Math., Vol. VIII pp. 181--189 Amer. Math. Soc., Providence,
R.I.



\bibitem {folkman} J. Folkman,
On the representation of integers as sums of distinct terms from a
fixed sequence, {\it Canad. J. Math.} 18 1966 643--655.

\bibitem{FF}
P. Frankl and Z.  F\"uredi, Solution of the Littlewood-Offord
problem in high dimensions. {\it Ann. of Math. } (2) 128 (1988), no.
2, 259--270.


\bibitem {freiman}  G. Freiman,
{ Foundations of a structural theory of set addition}. Translated
from the Russian. {\it Translations of Mathematical Monographs,} Vol
37. American Mathematical Society, Providence, R. I., 1973. vii+108
pp.

\bibitem{freiman2} G. Freiman, New analytical results in subset sum
problem, {\t Discrete mathematics} 114 (1993), 205-218.


\bibitem{GH} W. Gao and Y. O.  Hamidoune,
 On additive bases, {\it Acta Arith.} 88 (1999), no. 3, 233--237.

\bibitem{GHLS} W. Gao, Y.O. Hamidoune, A. Lladó and O. Serra,
Covering a finite abelian group by subset sums, {\it Combinatorica}
23 (2003), no. 4, 599--611.


\bibitem{Girko} V. Girko,  Circle law, (Russian)
 {\it Teor. Veroyatnost. i Primenen.}  29  (1984),  no. 4, 669--679.

  \bibitem{GT} F. G\"otze, A.N. Tikhomirov, On the circular law,
 \emph{preprint}.

 \bibitem{GT2} F. G\"otze, A.N. Tikhomirov, The Circular Law for Random
 Matrices,
 \emph{preprint}.



\bibitem{green} B. Green, The Cameron-Erdõs Conjecture
{\it Bull. London Math. Soc.} 36 (2004), no. 6, 769-778.


\bibitem{GR} B. Green and I. Ruzsa, Freiman's theorem in an arbitrary abelian
group,
{\it Jour. London Math. Soc.} 75 (2007), no. 1, 163-175.





\bibitem{GLOS}  J.  Griggs, J.  Lagarias, A. Odlyzko and J. Shearer,
 On the tightest packing of sums of vectors, {\it European J. Combin. }
 {4}   (1983), no. 3, 231--236.




\bibitem{HZ} Y. Hamidoune and G. Z\'emor, {
On zero-free subset sums}, {\it Acta Arithmetica} 78 (1996) no. 2,
143--152.




\bibitem{Kat} G. Katona,  On a conjecture of Erdös and a stronger form of Sperner's theorem.
{\it Studia Sci. Math. Hungar }{ \bf 1} 1966 59--63.

\bibitem{Kle} D. Kleitman,  On a lemma of Littlewood and Offord on the distributions of linear combinations of vectors,
{\it Advances in Math. }{\bf  5 } 1970 155--157 (1970).







\bibitem{Lipkin} {E. Lipkin} { On representation of $r-$powers by subset
sums}, {\it Acta Arithmetica}  52 (1989), 114-130.





\bibitem{LO} J. E.  Littlewood and A. C.  Offord,
On the number of real roots of a random algebraic equation. III.
{\it Rec. Math. [Mat. Sbornik] N.S. } {\bf 12} ,  (1943). 277--286.



\bibitem{EGZ} { P. Erd\H{o}s,
A. Ginzburg} and {A. Ziv}, { Theorem in the additive number theory},
{\it Bull. Res. Council Israel} 10F (1961), 41-43.

\bibitem{GPT} { W. D. Gao, A. Panigrahi} and {R. Thangadurai},
{\it On the structure of $p$-zero-sum-free sequences and its
application to a variant of Erd\H{os}-Ginzburg-Ziv theorem}. Proc.
Indian Acad. Sci. Vol. 115, No. 1 (2005), 67-77.


\bibitem{Hal} G. Hal\'asz,
Estimates for the concentration function of combinatorial number
theory and probability, {\it Period. Math. Hungar.} {\bf 8} (1977),
no. 3-4, 197--211.

\bibitem {Heg} N. Hegyv\'ari,
On the representation of integers as sums of distinct terms from a
fixed set, {\it Acta Arith. } 92 (2000), no. 2, 99--104.



\bibitem {LS} T. \L uczak and T. Schoen,
On the maximal density of sum-free sets, {\it Acta Arith.}  95
(2000), no. 3, 225--229.

\bibitem{MW} H. B. Mann and Y. F. Wou,
An addition theorem for the elementary abelian group of type
$(p,p)$, {\it Monatsh. Math. } 102 (1986), no. 4, 273--308.

\bibitem{Nathansonbook} M. Nathanson, Elementary methods in number
theory, Springer  2000.


\bibitem{NSV} H. H. Nguyen, E. Szemer\'edi and V.  Vu, Subset sums in
$\BZ_p$, {\it to appear in Acta Arithmetica}.

\bibitem{NV1} H. H. Nguyen and V. Vu, Classification theorems for sumsets modulo a prime  {\it submitted}.

\bibitem{NV2} H. Nguyen and V. Vu, On square-sum-free sets, {\it in
preparation}.

\bibitem{O1} J. E. Olson, Sums of sets of group elements,
{\it  Acta Arith.}  28 (1975/76), no. 2, 147--156.

\bibitem{O2} J. E. Olson,
An addition theorem modulo $p$, {\it J. Combinatorial Theory} 5 1968
45--52.



\bibitem{PZ} G. Pan and W. Zhou,
Circular law, Extreme singular values and potential theory, {\it
preprint}.

\bibitem{sar1} A. S\'ark\"ozi, Finite addition theorems I, {\it J. Number Theory,} 32, 1989, 114--130.

\bibitem{Sarkozy} {A.
S\'ark\"ozy}, { Finite Addition Theorems, II}, {\it Journal of
Number Theory}, 48 (1994), 197-218.

\bibitem{Sarhandbook} A. S\'ark\"ozy and C. Pomerance, Combinatorial
number theory, Chapter 20, Handook of Combinatorics (eds. R. Graham,
M. Gr\"otschel and  L. Lov\'asz), North-Holland 1995.



\bibitem{SS} A. S\'ark\"ozy and E. Szemer\'edi, \"Uber ein Problem von Erd\H{o}s und Moser,
{\it Acta Arithmetica}, 11 (1965) 205-208.




\bibitem{Stan} R. Stanley,  Weyl groups, the hard Lefschetz theorem, and the Sperner property, {\it SIAM J.
 Algebraic Discrete Methods} {\bf   1  } (1980), no. 2, 168--184.

\bibitem{SzEH} E. Szemer\'edi,
 On a conjecture of Erd\H os and Heilbronn, {\it Acta Arith.} 17 (1970) 227--229.



\bibitem{SzV1} { Endre Szemer\'edi} and { V. Vu },
{ Long arithmetic progression in sumsets and the number of x-free
sets}, {\it Proceeding of London Math Society}, 90(2005) 273-296.

\bibitem{SzV2} { E. Szemer\'edi} and { V.  Vu },
Finite and infinite arithmetic progression in  sumsets, {\it Annals
of Math,} 163 (2006), 1-35.



\bibitem{SzV3} { E. Szemer\'edi} and { V.  Vu },  Long
arithmetic progressions in sumsets: Thresholds and Bounds, {\it
Journal of the A.M.S,} 19 (2006), no 1, 119-169.







\bibitem{TVbook} T. Tao and V.    Vu, Additive Combinatorics, Cambridge Univ. Press, 2006.
\bibitem{TVcir} T. Tao and V. Vu, Random matrices: The Circular
 Law, {\it Communication in Contemporary Mathematics} 10 (2008), 261-307.

\bibitem{TVinverse} T. Tao and V. Vu, Inverse Littlewood-Offord theorems
and the condition number of random matrices, {\it to appear in
Annals of Mathematics}.

\bibitem{TVinverse08} T. Tao and V. Vu, {\it paper in preparation}.


\bibitem{Vudir} V. Vu,   Structure of large incomplete sets in abelian
groups, {\it to appear in Combinatorica}.


\bibitem{Wig} Wigner, On the distribution of the roots of certain symmetric matrices,
  {\it Annals of Mathematics} (2)  67  1958 325--327.


\end{thebibliography}
\end{document}